\documentclass[12pt,a4paper,reqno]{amsart}
\newtheorem{theorem}{Theorem}
\newtheorem{lemma}{Lemma}
\newtheorem{corollary}{Corollary}
\begin{document}

\title{Turan Extremum Problem for Periodic Function with Small
Support}
\author{D.V.~Gorbachev, A.S.~Manoshina}
\begin{abstract}
We consider an extremum problem posed by Turan. The aim of this problem is to find
a maximum mean value of 1-periodic continuous even function such that sum of
Fourier coefficient modules for this function is equal to 1 and support of this
function lies in $[-h,h]$, $0<h\le 1/2$. We show that this extremum problem for rational
$h=p/q$ is equivalent two finite-dimensional linear programming problems. Here
there are exact results for rational $h=2/q$, $h=p/(2p+1)$, $h=3/q$, and asymptotic equalities.
\end{abstract}
\date{11/27/2001}
\thanks{This research was supported by the RFBR unedr grant no. 00-01-00644,
01-01-06190 and 02-01-06563.
\endgraf
This work was published in Tchebyshev Collection, TPSU (Russia), \textbf{2} (2001), 31--40.
}
\email{dvg@uic.tula.ru; manoshin@tula.net}
\urladdr{http://home.uic.tula.ru/$\sim$gd030473}
\maketitle

Suppose $0<h\le1/2$; then by $K(h)$ we denote a class of continuous
1-periodic even functions $f(x)$, where $f(x)$ satisfies following conditions:

1) $f(x)=\sum_{n=0}^{\infty}a_n\cos(2\pi nx)$;

2) $a_n\ge0$ for all $n=0,1,2,\dots$;

3) $f(0)=\sum_{n=0}^{\infty}a_n=1$;

4) $f(x)=0$, for $h\le|x|\le1/2$.

For arbitrary $h$ this class $K(h)$ is non-empty. For example, if 1-periodic
function $\varphi(x)$ possesses the value $\max\{1-h^{-1}|x|,0\}$ on the interval
$[-1/2,1/2]$ and has expansion in Fourier series \cite{Stechkin}
$$\varphi_h=h+2h\sum_{n=1}^{\infty}\left(\frac{\sin(\pi nh)}{\pi
nh}\right)^{2}\cos(2\pi nx),$$ then this function belongs to $K(h)$.

We estimate the following value $$A(h)=\sup_{f\in K(h)}a_0=\sup_{f\in
K(h)}\int_{-h}^hf(x)dx.$$

In 1970 P. Turan posed this extremum problem in connection with applications in the
theory of number. This problem is connection with van der Corput sets. Other
applications in the analytic number theory your can see in \cite{Konyagin}.

In \cite{Stechkin}, S.B. Stechkin  proved that $A(1/q)=1/q$ for $h=1/q$
($q=2,3,\dots$). In this case the extremal function is $\varphi(x)$.

For $h\ne1/q$ $A(h)>h$ (A.Y. Popov, verbal information).

\begin{theorem}[\cite{Gorbachev}]
$A(h)=h+O(h^3)$ as $h\to 0$. Asymptotically an extremal function is $\varphi(x)$.
The constant in $O$ is no more than $16$.
\end{theorem}


Bellow we consider the problem of Turan for arbitrary ration $h$. Our main results
are following.

Suppose $p,q\in \mathbb{N}$, $(p,q)=1$ and $2p\le q$. We consider two problems of
linear programming.
\begin{align}\text{1. If}\quad &\sum_{r=0}^{q-1}s_r=1, \notag
\\ &\sum_{r=0}^{q-1}s_r\cos(2\pi rk/q)=0,\quad k=p,p+1,\dots,q-p, \notag \\
&s_r\ge 0,\quad r=0,1,\dots,q-1; \notag \\ &S_1(p,q)=\max_{s}s_0\quad\text{is to
be found.} \notag \\ \notag \\ \text{2. If}\quad&1+2\sum_{k=1}^{p-1}b_k\cos(2\pi
rk/q)\ge0, \quad r=0,1,\dots q-1; \notag
\\ &S_2(p,q)=\max_{b}\frac{1}{q}\left(1+2\sum_{k=1}^{p-1}b_k\right)\quad \text{is
to be found.} \notag
\end{align}

\begin{theorem} \label{theorem2}
Suppose $p,q\in \mathbb{N}$, $(p,q)=1$ and $2p\le q$; then
$$A(p/q)=S_1(p,q)=S_2(p,q).$$ Extremal function is $1$-periodic piece linear function
$\varphi_{p,q}(x)$ with vertex at points $$(\pm k/q,b_k^{*}),\quad
k=0,1,\dots,p\quad b_0^{*}=1,\quad b_p^{*}=0.$$ and $$\varphi_{p,q}(x)=
\begin{cases}
b_k^{*}+(b_k^{*}-b_{k+1}^{*})(k-q|x|), & k/q\le|x|<(k+1)/q,\\ 0, & p/q\le|x|<1/2.
\end{cases}$$
Here $b_k^{*}$ $(k=1,2,\dots,p-1)$ are values of the optimal schedule in linear
programming problem~$2$.
\end{theorem}

To prove Theorem \ref{theorem2}, we need several lemmas.

\begin{lemma}
$A(p/q)\le S_1(p,q)$.
\end{lemma}
\emph{Proof.} By $f(x)$ denote an arbitrary function of $K(p/q)$. For $f(x)$
following conditions are hold.

1. $f(k/q)=0$, $k=p,p+1,\dots,q-p$ (from fourth property of $K(p/q)$). For expansion in
Fourier series we get $$\sum_{n=0}^{\infty}a_n\cos(2\pi nk/q)=0, \quad
k=p,p+1,\dots,q-p.\eqno{(1)}$$

Now we have $n=q\nu+r$, where $\nu=0,1,2,\dots$, $r=0,1,\dots,q-1$ are quotient and
residue of division $n/q$ respectively. Then $\cos(2\pi nk/q)=\cos(2\pi rk/q)$.
Using (1), we get $$\sum_{r=0}^{q-1}s_r\cos(2\pi rk/q), \quad
k=p,p+1,\dots,q-p,\eqno{(2)}$$ where $s_r=\sum_{\nu=0}^{\infty}a_{q\nu +r}$. By
second property of $K(h)$, it follows that $$s_r\ge 0,\quad
r=0,1\dots,q-1.\eqno{(3)}$$

2. Since $\sum_{n=0}^{\infty}a_n=1$, we obtain $$\sum_{r=0}^{q-1}s_r=1.\eqno{(4)}$$

3. Clearly, $a_0\le s_0$.

Since (2),(3),(4) are conditions of linear programming problem 1, we see that
$s_0\le S_1(p,q)$. Then $A(p/q)\le S_1(p,q)$.

\begin{lemma}
$A(p/q)\ge S_2(p,q)$.
\end{lemma}
\emph{Proof.} By $K_0(p/q)$ denote subclass of $K(p/q)$ such that it consists of
piece linear functions with vertexes at points $$\left(\pm k/q,b_k\right),\quad
k=0,1,\dots,p,\quad b_0=1,\quad b_p=0.$$ Any such function is given by
$$\varphi(x)=
\begin{cases}
b_k+(b_k-b_{k+1})(k-q|x|), & k/q\le|x|<(k+1)/q,\\ 0, & p/q\le|x|<1/2.
\end{cases}\eqno{(5)}$$
Its expansion in Fourier series is
$$\varphi(x)=\sum_{n=0}^{\infty}\alpha_{n}\cos(2\pi nx),$$ where
$$\alpha_{0}=\frac{1}{q}\left(1+2\sum_{k=1}^{p-1}b_k\right),\quad
\alpha_{n}=2q\left(\frac{\sin(\pi n/q)}{\pi
n}\right)^2\left(1+2\sum_{k=1}^{p-1}b_k\cos(2\pi nk/q)\right).$$ For $\varphi(x)\in
K_0(p/q)$ it is necessary to have $\alpha_{n}\ge 0,$ ($n=0,1,2,\dots$), i.e.
$$1+2\sum_{k=1}^{p-1}b_k\cos(2\pi nk/q)\ge 0,\quad n=0,1,2,\dots$$ or
$$1+2\sum_{k=1}^{p-1}b_k\cos(2\pi rk/q)\ge 0,\quad r=0,1,\dots q-1.$$

The subclass $K_0(p/q)$ is non-empty. In particular $\varphi_{p/q}(x)\in K_0(p/q)$.
We shall seek function from $K_0(p/q)$ such that its coefficient $\alpha_0$ is
highest possible. For that we solve
$$\alpha_0=\frac{1}{q}\left(1+2\sum_{k=1}^{p-1}b_k\right)\to \max$$ under
conditions $$1+2\sum_{k=1}^{p-1}b_k\cos(2\pi rk/q)\ge 0,\quad r=0,1,\dots q-1.$$
This problem is the linear programming problem 2.

Let $b_k^{*}$ (k=1,2,\dots,p-1) be an optimal schedule in the problem 2. We define
function $\varphi_{p,q}(x)$ from $K_0(p/q)\subset K(p/q)$ by formula (5) with
coefficients $b_k^{*}$. Its mean value is $S_2(p,q)$. Since $\varphi_{p,q}\in
K(p/q)$ we have $A(p/q)\ge S_2(p,q)$.

\begin{lemma}
$S_1(p,q)=S_2(p,q)$.
\end{lemma}
\emph{Proof.} Let us consider trigonometric sums from problem 2
$$b_k=\sum_{r=0}^{q-1}s_r\cos(2\pi rk/q),\quad k=0,1,\dots q-1,$$ where
$$b_0=1,\quad b_p=b_{p+1}=\dots=b_{q-p}=0.\eqno{(6)}$$ Using discrete
cosine Fourier inversion, we get
$$\frac{s_r+s_{q-r}}{2}=\frac{1}{q}\sum_{k=0}^{q-1}b_k\cos(2\pi rk/q),\quad
r=0,1,\dots,q-1,\quad s_0=s_q.$$ Since $\cos(2\pi r(q-k)/q)=\cos(2\pi rk/q)$, by
(6) it follows that
$$\frac{s_r+s_{q-r}}{2}=\frac{1}{q}\left(1+2\sum_{k=1}^{p-1}b_k\cos(2\pi
rk/q)\right).$$ Therefore the finding value $S_1(p,q)$ is equivalent following
problem: $$\max s_0=\max_{b}\frac{1}{q}\left(1+2\sum_{k=1}^{p-1}b_k\right)$$ under conditions
$$q\frac{s_r+s_{q-r}}{2}=1+2\sum_{k=1}^{p-1}b_k\cos(2\pi rk/q)\ge 0,\quad r=0,1,\dots
q-1.$$ But this is the linear programming problem 2 and lemma follows.

Thus from lemmas 1--3, we get Theorem 2. Using this theorem, we can find value
$A(p/q)$ for fast values $p$ and $q$.

\begin{theorem}
For $p=2$ and $h=2/q$, $q=5,7,\dots$ $$A(2/q)=\frac{1+\cos(\pi/q)}{q\cos(\pi/q)}.$$
\end{theorem}
\emph{Proof.} By Theorem 2, so that $A(2/q)=S_2(2,q)$. We have
$$A(2/q)=\frac{1}{q}\max_{b}(1+2b),\quad 1+2b\cos(2\pi r/q)\ge0, r=0,1\dots,q-1.$$
For $r=0,1,\dots,q-1$ and odd $q$ we have $\min(\cos(2\pi r/q))=-\cos(\pi/q)$.
Therefore $$\max 2b=2b_1^{*}=\frac{1}{\cos(\pi/q)}\eqno{(7)}$$ and
$$A(2/q)=\frac{1}{q}\left(1+\frac{1}{\cos(\pi/q)}\right)=\frac{1+\cos(\pi/q)}{q\cos(\pi/q)}.$$

An extremal function is the piece linear function $\varphi_{2,q}(x)$ (5) with
vertex at points ($0,1$), ($\pm 1/q,b_1^{*}$), ($\pm 2/q,0$) (7).

\begin{theorem}
For $h=\frac{p}{2p+1}$, $p=1,2,\dots$
$$A\left(\frac{p}{2p+1}\right)=\frac{\cos(\frac{\pi}{2p+1})}{1+\cos(\frac{\pi}{2p+1})}.$$
\end{theorem}
\emph{Proof.}
From Theorem 2, we get $A(\frac{p}{2p+1})=S_1(p,2p+1)$ and
$$
A\left(\frac{p}{2p+1}\right)=\max_{s}s_0
$$
under conditions
\begin{align}
&\sum_{r=0}^{q-1}s_{r}=1, \tag{8} \\
&\sum_{r=0}^{q-1}s_{r}\cos(2\pi rk/q)=0,\quad k=p,p+1, \tag{9}\\
&s_{r}\ge 0,\quad r=0,1,\dots q-1.\tag{10}
\end{align}
From (8) we have $s_{1}=1-s_{0}-\sum_{r=2}^{q-1}s_{r}$. For $k=p$ and (9)
$$
s_{0}+\left(1-s_{0}-\sum_{r=2}^{q-1}s_{r}\right)\cos(2\pi
p/q)+\sum_{r=2}^{q-1}s_{r}\cos(2\pi rp/q)=0.
$$
Therefore
\begin{equation}
s_{0}\left(1-\cos(2\pi p/q)\right)=-\cos(2\pi
p/q)+\sum_{r=2}^{q-1}s_r\left(\cos(2\pi p/q)-\cos(2\pi rp/q)\right).\tag{11}
\end{equation}
Since for $q=2p+1$ $$\cos\left(\frac{2\pi
p}{q}\right)=-\cos\left(\frac{\pi}{2p+1}\right)$$
and
\begin{align}
\cos\left(\frac{2\pi p}{q}\right)-\cos\left(\frac{2\pi rp}{q}\right)=
-\cos\left(\frac{\pi}{2p+1}\right)-(-1)^r\cos\left(\frac{\pi r}{2p+1}\right)\le0,
\notag\\
r=2,3,\dots, q-1,\notag
\end{align}
from (10)--(11) we get
\begin{align}
&s_0\left(1+\cos\left(\frac{\pi}{2p+1}\right)\right)=
\cos\left(\frac{\pi}{2p+1}\right)+\notag\\
&\qquad+\sum_{r=2}^{q-1}s_r\left(-\cos\left(\frac{\pi}{2p+1}\right)-
(-1)^r\cos\left(\frac{\pi r}{2p+1}\right)\right)
\le\cos\left(\frac{\pi}{2p+1}\right)\tag{12}
\end{align}
Therefore $$\max_s
s_0\le\frac{\cos\left(\frac{\pi}{2p+1}\right)}{1+\cos\left(\frac{\pi}{2p+1}\right)}.
$$
In (12) we have equality for
\begin{align}
s_0^{*}&=\frac{\cos\left(\frac{\pi}{2p+1}\right)}{1+\cos\left(\frac{\pi}{2p+1}\right)}>0,\notag\\
s_1^{*}&=\frac{1}{1+\cos\left(\frac{\pi}{2p+1}\right)}>0,\notag\\
s_r^{*}&=0,\quad r=2,3,\dots,q-1.\notag
\end{align}
Thus $$A\left(\frac{p}{2p+1}\right)=\frac{\cos(\frac{\pi}{2p+1})}{1+\cos(\frac{\pi}{2p+1})}$$
and extremal function is piece linear function
$\varphi_{p,2p+1}(x)$ with vertex at points $(\pm k/q,b_k^{*})$, $k=0,1,\dots,p,$
where (see proof of lemma 3)
$$b_k^{*}=\frac{\cos\left(\frac{\pi}{2p+1}\right)+\cos\left(\frac{2\pi k}{2p+1}\right)}{1+\cos\left(\frac{\pi}{2p+1}\right)}.$$

Notice that asymptotically (as $p\to \infty$) an extremal function is $\cos^2(\pi x)$.

\begin{theorem}
For $p=3$ and $h=3/q$, $q=7,8,10,11,\dots$
$$A(3/q)=\frac{1}{q}\left(1+\frac{1-2(\cos(2\pi
r_0/q)+\cos(2\pi(r_0+1)/q))}{1+2\cos(2\pi r_0/q)\cos(2\pi(r_0+1)/q)}\right),$$
where $r_0=[q/3]$ is an integer part of $q/3$.
\end{theorem}
\emph{Proof.} From Theorem 2, we get
$$A(3/q)=S_2(3,q)=\frac{1}{q}\max_{\alpha}\alpha(0),\quad \alpha(r)\ge 0,\quad
r=0,1,\dots,q-1,\eqno{(13)}$$ where $\alpha(r)$ is trigonometric polynomial
$$\alpha(r)=1+2b_1\cos(2\pi r/q)+2b_2\cos(4\pi r/q). \eqno{(14)}$$

We define a quadrature formula
$$1=\frac{1}{q}\sum_{r=0}^{q-1}\alpha(r)=\gamma_0\alpha(0)+\gamma_1\alpha(r_0)+\gamma_2\alpha(r_0+1),\quad
r_0=[q/3].\eqno{(15)}$$

For equality (15) for any trigonometric polynomial (14), it is sufficient to have
(15) for basis functions $1$, $\cos(2\pi r/q)$, $\cos(4\pi r/q)$:
\begin{align}1=&\frac{1}{q}\sum_{r=0}^{q-1}1=\gamma_0+\gamma_1+\gamma_2, \notag
\\ 0=&\frac{1}{q}\sum_{r=0}^{q-1}\cos(2\pi r/q)=\gamma_0+\gamma_1\cos(2\pi r_0/q)+\gamma_2\cos(2\pi(r_0+1)/q), \notag \\
0=&\frac{1}{q}\sum_{r=0}^{q-1}\cos(4\pi r/q)=\gamma_0+\gamma_1\cos(4\pi
r_0/q)+\gamma_2\cos(4\pi(r_0+1)/q). \notag
\end{align}

Since $\gamma_0\ne0$, we get combined equations $$\begin{cases} \gamma_1'\cos(2\pi
r_0/q )+\gamma_2'\cos(2\pi(r_0+1)/q)=-1,\\\gamma_1'\cos(4\pi
r_0/q)+\gamma_2'\cos(4\pi(r_0+1)/q)=-1,\end{cases}$$ where
$\gamma_1'=\gamma_1/\gamma_0$ and $\gamma_2'=\gamma_2/\gamma_0.$

We have $$\gamma_i'=\frac{\Delta_i}{\Delta},\quad i=1,2.\eqno{(16)}$$ Here
\begin{align}\Delta=&\cos(2\pi r_0/q)\cos(4\pi(r_0+1)/q)-\cos(2\pi(r_0+1)/q)\cos(4\pi r_0/q)= \notag \\
=&(\cos(2\pi(r_0+1)/q)-\cos(2\pi r_0/q))(1+2\cos(2\pi r_0/q)\cos(2\pi(r_0+1)/q))<0
\notag
\end{align}
(from inequalities $q/3-1<r_0<r_0+1<q/3+1$),
\begin{align}\Delta_1=&\cos(2\pi(r_0+1)/q)-\cos(4\pi(r_0+1)/q)= \notag \\
=&2\sin(\pi(r_0+1)/q)\sin(3\pi(r_0+1)/q))<0 \notag
\end{align}
(from inequalities $q/3<r_0+1<q/3+1$),
\begin{align}\Delta_2=&\cos(4\pi r_0/q)-\cos(2\pi r_0/q)= \notag \\
=&-2\sin(\pi r_0/q)\sin(3\pi r_0/q))<0 \notag
\end{align}
(from inequalities $q/3-1<r_0<q/3$).

Thus, $\gamma_1',\gamma_2'>0$. Therefore all weights $\gamma_i$ ($i=0,1,2$) have
like signs. From $\gamma_0+\gamma_1+\gamma_2=1$ we get $\gamma_i>0$ and
$\gamma_1'+\gamma_2'=(1-\gamma_0)/\gamma_0$. From (16) we have
$$\frac{1}{\gamma_0}=1+\gamma_1'+\gamma_2'=1+\frac{1-2(\cos(2\pi
r_0/q)+\cos(2\pi(r_0+1)/q))}{1+2\cos(2\pi r_0/q)\cos(2\pi (r_0+1)/q)}.\eqno{(17)}$$

Using the quadrature formula (15) for arbitrary trigonometric polynomial (14)
($\alpha(r)\ge 0$, $r=0,1,\dots,q-1$), we get
$$1=\frac{1}{q}\sum_{r=0}^{q-1}\alpha(r)=\gamma_0\alpha(0)+\gamma_1\alpha(r_0)+\gamma_2\alpha(r_0+1).$$

Since $\gamma_1\alpha(r_0)\ge 0$ and $\gamma_2\alpha(r_0+1)\ge 0$, it follows that
$\gamma_0\alpha(0)\le 1$ and $$\alpha(0)\le\frac{1}{\gamma_0}.$$

The equality reaches for $\alpha^{*}(r)=1+2b_1^{*}\cos(2\pi r/q)+2b_2^{*}\cos(4\pi
r/q),$ where $$b_1^{*}=-\frac{\cos(2\pi r_0/q)+\cos(2\pi(r_0+1)/q)}{1+2\cos(2\pi
r_0/q)\cos(2\pi(r_0+1)/q)}>0,$$ $$b_2^{*}=\frac{1/2}{1+2\cos(2\pi
r_0/q)\cos(2\pi(r_0+1/q))}>0. \eqno{(18)}$$

An extremal function is the piece linear function $\varphi_{3,q}(x)$ (5) with
vertex at points ($0,1$), ($\pm 1/q,b_1^{*}$), ($\pm 2/q,b_2^{*}$), ($\pm 3/q,0$)
(18). Theorem is proved.

Notice that the main idea of using quadrature formula (15) belong to V.I.~Ivanov.
In \cite{Manoshina} there is an other proof.

Finally, we obtain asymptotic equalities.

\begin{corollary}
For $h=2/q\to 0$, $q=2\nu+1$, $\nu=2,3,\dots$
$$A(h)=h+\frac{\pi^2}{16}h^3+\frac{5\pi^4}{768}h^5+O(h^7).$$ For $h=3/q\to 0$,
$q=3\nu+1$, $\nu=2,3,\dots$
$$A(h)=h+\frac{16\pi^2}{243}h^3-\frac{16\pi^3\sqrt{3}}{2187}h^4+\frac{448\pi^4}{59049}h^5+O(h^6),$$
for $h=3/q\to 0$, $q=3\nu+2$, $\nu=2,3,\dots$
$$A(h)=h+\frac{16\pi^2}{243}h^3+\frac{16\pi^3\sqrt{3}}{2187}h^4+\frac{448\pi^4}{59049}h^5+O(h^6).$$
\end{corollary}

Therefore in Theorem 1 we can not replace O-large by o-small. We express one's
thanks to V.I.~Ivanov for useful discussions.

\end{document}